\documentclass[onecolumn]{revtex4}  %% REVTeX 4.0

\usepackage{amsmath,amsfonts,amssymb}
\usepackage{graphicx}

\graphicspath{{.}{figures/}}

\newcommand\eqnref[1]{(\ref{#1})}
\newcommand\figref[1]{Fig.~\ref{#1}}

\newcommand{\iu}   {\mathrm{i}}     %  Imaginary Unit
\newcommand{\Deltarm}   {\mathrm{\Delta}}

\newcommand{\Omegarm}   {\mathrm{\Omega}}

\newcommand{\bfk}   {\mathbf{k}}

\newcommand{\bfr}   {\mathbf{r}}

\begin{document}

\title{A ``Trefftz Machine'' for Absorbing Boundary Conditions}

\author{Igor Tsukerman}
\affiliation{Department of Electrical and Computer Engineering, The
University of Akron, OH 44325-3904, USA}

\email{igor@uakron.edu}

\begin{abstract}
The paper presents an automatic generator of approximate nonreflecting boundary conditions,
analytical and numerical, for scalar wave equations. This generator has two main ingredients. 
The first one is a set of local Trefftz functions --  outgoing waves
approximating the solution in the vicinity of a given point of the exterior boundary
of the computational domain.
The second ingredient is a set of linear test functionals (degrees of freedom). 
One example of such functionals is the nodal values of the solution at a set of grid points; 
in that case, one obtains a numerical condition -- a finite difference scheme at the boundary.
Alternatively, the functionals may involve derivatives or integrals
of the solution, in which case the proposed ``Trefftz machine'' yields analytical
nonreflecting conditions. Corners and edges are treated algorithmically
the same way as straight boundaries. With specific choices of bases and degrees of freedom,
the machine produces classical conditions such as Engquist-Majda and Bayliss-Turkel.
For other choices, one obtains a variety of analytical and numerical conditions,
a few of which are presented as illustrative examples.
\end{abstract}

%\begin{keywords}
%\end{keywords}

%\date{Received: XXX}

% \today

\maketitle

%%%%%%%%%%%%%%%%%%%%%%%%%%%%%%%%%%%%%%%%%
\section{Introduction}
\label{sec:introduction}
%%%%%%%%%%%%%%%%%%%%%%%%%%%%%%%%%%%%%%%%%
%

The critical role of artificial boundary conditions for finite difference or finite element solution
of wave problems is well recognized. 
This subject is four decades old, is vast and includes various types of Perfectly Matched
Layers (PML, see e.g. \cite{Berenger96,Teixeira98,Sacks95,Gedney96,Collino98,Becache02,Becache04}) 
and absorbing boundary conditions
(e.g. \cite{Engquist77,Higdon86,Higdon87,Bayliss-Turkel80,Bayliss-Gunzburger-Turkel82,Hagstrom-Hariharan98,
Givoli01,Givoli-Neta03,Hagstrom03,Hagstrom-Wartburton04,Hagstrom08,Hagstrom-Warburton-Givoli10,
Zarmi-Turkel13}).
A number of excellent reviews are available, e.g. \cite{Givoli04,Tsynkov98,Hagstrom99},
so here I highlight only two classical ideas directly related to the material of this paper.

We shall consider the scalar wave equation either in the frequency domain
\begin{equation}\label{eqn:del2-u-plus-k2u}
   \nabla^2 u + k_0^2 u \,=\, f  ~~~\mathrm{in}~~R^n, ~~ n=1,2,3;
   ~~ \mathrm{supp} \, f \subset \Omegarm \subset R^n
\end{equation}
or, alternatively, in the time domain
\begin{equation}\label{eqn:del2-u-d2u-dt2}
    v^2 \, \nabla^2 u - \partial^2_{tt} u ~=~ f   ~~~\mathrm{in}~~R^n, ~~ n=1,2,3;
   ~~ \mathrm{supp} \, f \subset \Omegarm \times [0, \infty)
\end{equation}
As indicated in these equations, sources $f$ are assumed to be confined to 
a bounded domain $\Omegarm$ in space. In \eqnref{eqn:del2-u-plus-k2u}, $k_0$  is a given positive parameter 
(the wavenumber). In \eqnref{eqn:del2-u-d2u-dt2}, $v$ is the velocity of waves,
for simplicity assumed to be position-independent, although the approach of this paper 
can be extended to more complex cases. When convenient for
analysis, $v$ will be normalized to unity. We shall deal primarily with 2D problems, 
although all ideas can be extended to 3D.
Let us assume that $\partial \Omegarm$ is a rectangle (a parallelepiped in 3D); conditions at the corners
(and edges) will \emph{not} be ignored.

Problem \eqnref{eqn:del2-u-plus-k2u} requires radiation boundary conditions 
(e.g. Sommerfeld) at infinity, but our task is 
to replace these theoretical conditions with approximate but accurate and practical ones
on the exterior surface $\partial \Omegarm$ away from the sources.
This is to be done in such a way that the solution subject to these artificial conditions 
be by some measure close to the true solution in $\Omegarm$.
Similarly, we shall seek approximate boundary conditions 
on $\partial \Omegarm$ for problem \eqnref{eqn:del2-u-d2u-dt2} as well.
Initial conditions for \eqnref{eqn:del2-u-d2u-dt2} are assumed to be given
and are tangential to our analysis.

Let us first consider a straight artificial boundary in 2D, 
for convenience at $x = 0$, with the computational domain $\Omegarm$ 
situated on the positive $x$ side.
One classical nonreflecting condition, due to Engquist \& Majda \cite{Engquist77}, 
follows from the dispersion relation
\begin{equation}\label{eqn:dispersion-relation-kx-plus-sqrt}
    k_x^2 + k_y^2 - k^2 = 0 ~~~~ \Leftrightarrow ~~~~ k_x = - k \, \sqrt{1 - k_y^2/k^2}
\end{equation}
which holds for problem \eqnref{eqn:del2-u-plus-k2u} with $k=k_0$
and for problem \eqnref{eqn:del2-u-d2u-dt2} with $k = \omega / v$, $\omega \in (-\infty, \infty)$.
The negative sign of $k_x$ in \eqnref{eqn:dispersion-relation-kx-plus-sqrt} corresponds 
to outgoing waves (waves moving in the $-x$ direction) under the $\exp(-\iu \omega t)$ phasor convention.

If instead of the square root \eqnref{eqn:dispersion-relation-kx-plus-sqrt}
contained a rational fraction of $k_x, k_y, k_0$, then the inverse transform 
of the corresponding dispersion relation would be an exact nonreflecting 
boundary condition involving a combination of $x$, $y$ and $t$ derivatives.
It is then clear that a sequence of approximate absorbing conditions
can be derived using Taylor or Pad\'{e} approximations
of the square root and inverse-transforming these relationships back to real space;
see e.g. \cite{Cai-book13} for details.
(Engquist \& Majda's analysis is ultimately equivalent but cast 
in the language of pseudodifferential operators.)

Another classical idea, due to Bayliss \& Turkel \cite{Bayliss-Turkel80,Zarmi-Turkel13}, involves
a cylindrical (2D) or spherical (3D) harmonic expansion of radiated fields.
A sequence of differential operators annihilating progressively higher numbers of terms
in this expansion constitutes absorbing conditions of progressively higher orders.

One well recognized shortcoming of these classical methods is their reliance on high-order
derivatives that are difficult to deal with in numerical simulations.
To overcome this deficiency, several clever reformulations
have been proposed \cite{Givoli-Neta03,Hagstrom-Hariharan98,Hagstrom-Wartburton04,
Hagstrom-Warburton-Givoli10}, with a sequence of auxiliary variables
on the exterior boundary instead of high-order derivatives.
Methods of this type will remain out of the scope of the present paper.
Rather, its goal is to devise a new ``machine'' for generating approximate absorbing schemes
that include, but are certainly not limited to,
the classical  Engquist-Majda and Bayliss-Turkel conditions.
Several examples of such schemes are presented in subsequent sections.

The ``machine'' has two main ingredients. The first one is a set of local basis functions
$\psi_\alpha$ ($\alpha = 1,2,\ldots, n$) approximating the
solution near a given point on the exterior boundary. These functions
can be chosen as outgoing waves of the form $g(\hat{\bfk} \cdot \bfr - t)$,
where $g$ is a given function (e.g. sinusoidal) and $\hat{\bfk}$
is a unit vector at an acute angle to the outward normal on $\partial \Omegarm$.
The second ingredient is a set of $m$ degrees of freedom (dof) -- linear functionals
$l_\beta (u)$ ($\beta = 1,2,\ldots, m$); $m$ is \textit{not} usually equal to $n$.

To elaborate, let the exact solution be approximated locally as a linear combination
\begin{equation}\label{eqn:uh-eq-c-psi}
    u_a = \sum_\alpha c_\alpha \psi_\alpha = \underline{c}^T \underline{\psi}
\end{equation}
where $\underline{c}$ is a Euclidean coefficient vector and $\underline{\psi}$
is a vector of basis functions. (Vectors are underlined to distinguish
them from other entities.) The functions and coefficients can be real-valued
or complex-valued, as will be clear from the context.
Coefficients $\underline{c}$ may be different at different boundary points,
but for simplicity of notation this is not explicitly indicated.

We are looking for a suitable boundary condition of the form
\begin{equation}\label{eqn:s-beta-l-beta-eq-0}
    \sum_\beta s_\beta l_\beta (u_a) ~=~ 0
\end{equation}
where $\underline{s} \in \mathbb{R}^m$ (or $\mathbb{C}^m$ in the complex case)
is a set of coefficients (``scheme'') to be determined.
We require that the scheme be exact for any $u_a$, i.e. for any linear combination of basis functions:
$$
    \sum_\beta s_\beta l_\beta \left( \sum_\alpha c_\alpha \psi_\alpha \right)  ~=~ 0
$$
or in matrix form
$$
    \underline{c}^T N^T \underline{s} ~=~ 0
$$
where $N^T$ is an $n \times m$ matrix with entries $N^T_{\alpha \beta} = l_\beta (\psi_\alpha)$.
Since the above equality is required to hold for all $\underline{c}$, one must have
\begin{equation}\label{eqn:s-in-null-Nt}
    \underline{s} \in \mathrm{Null} ~ N^T
\end{equation}

This whole development is completely analogous to that of FLAME
\cite{Tsukerman06,Tsukerman-book07,Tsukerman05},
where the goal is to construct a finite difference scheme rather than an absorbing condition.
The dof in FLAME are the nodal values of the solution on a given grid stencil 
\footnote{
(i) It is for the sake of compatibility of notation with the FLAME papers
that the matrix has been denoted with $N^T$ rather than just $N$.
(ii) In this paper, the term ``stencil'' means the set of nodes over which
a difference scheme is defined, not the coefficients of that scheme.}.
It is, however, interesting to bring more general linear functionals into consideration,
which is done in subsequent sections.

Gratkowski \cite{Gratkowski-book09} uses similar ideas to derive
analytical boundary conditions, albeit for static problems only 
and without the nullspace formula \eqnref{eqn:s-in-null-Nt}.
As multiple examples below and in \cite{Tsukerman-book07,Tsukerman05,Tsukerman06,
Tsukerman-PBG08,Tsukerman-JOPA09} demonstrate, this formula,
despite its simplicity, is rich and leads to a variety of useful schemes,
not only numerical as in the previous publications, but also analytical
as in the present paper.

One important measure of the quality of the boundary condition
is the reflection coefficient $R(\theta)$, defined as follows.
Consider an outgoing complex-exponential wave
$u_o = A_o \exp \left( \iu (-x \cos \theta_\alpha + y \sin \theta_\alpha + t) \right)$
and the corresponding reflected wave
$u_r = A_r \exp \left( \iu (x \cos \theta_\alpha + y \sin \theta_\alpha + t) \right)$,
where $A_o$, $A_r$ are complex amplitudes. Further, let the absorbing condition be defined 
by a set of coefficients $s_\beta$. Then, by definition, $R$ satisfies
$$
   \sum_\beta s_\beta l_\beta(u_o + Ru_r) \,=\, 0
$$
or
\begin{equation}\label{eqn:R-eq-lu-inc-over-lu-ref}
    R \,=\, -\frac{\sum_\beta s_\beta l_\beta(u_o)}
    {\sum_\beta s_\beta l_\beta(u_r)}
\end{equation}
%

%%%%%%%%%%%%%%%%%%%%%%%%%%%%%%%%%%%%%%%%%
\section{Example: Basis of Cylindrical Harmonics, Derivatives as dof}
\label{sec:Cyl-harmonic-basis}
%%%%%%%%%%%%%%%%%%%%%%%%%%%%%%%%%%%%%%%%%
%
Consider the 2D Helmholtz equation \eqnref{eqn:del2-u-plus-k2u}.
As we shall see in this section, the ``machine'' described above produces, 
with a natural choice of basis functions and degrees of freedom,
the classical Bayliss-Turkel conditions.
Indeed, the scattered field outside $\Omegarm$ can be expanded into cylindrical harmonics as
\begin{equation}\label{eqn:u-eq-sum-cn-hn}
   u(\bfr) \,=\, \sum_{n=-\infty}^{\infty} c_n h_{|n|}(k_0 r) \exp(\iu n \theta),
\end{equation}
where $h_n$ is the Hankel function (of the first kind, under the $\exp(-\iu \omega t)$
phasor convention for time-harmonic functions).
It is convenient to replace Hankel functions with their asymptotic expansions at infinity, viz.:
$$
   h_n(w) \,=\, \left( \frac{2}{\pi w} \right)^\frac12
   \exp \left(\iu \left(w - \frac{n \pi}{2} - \frac{\pi}{4} \right) \right)
   \sum_{l=0}^{\infty} \frac{a_l}{w^l}
$$
with some coefficients $a_l$, expressions for which are rather cumbersome and unimportant for our purposes.
Substituting this Hankel expansion into series \eqnref{eqn:u-eq-sum-cn-hn} for $u$,
one obtains
$$
   u(\bfr) \,=\, \left( \frac{2}{\pi k_0 r} \right)^\frac12
   \sum_{n=-\infty}^{\infty} c_n
   \exp \left(\iu \left(k_0 r - \frac{|n| \pi}{2} - \frac{\pi}{4} \right) \right)
   \exp(\iu n \theta) \, \sum_{l=0}^{\infty} \frac{a_l}{(k_0 r)^l} 
$$
\begin{equation}\label{eqn:u-eq-asymptotic-series}
  \sim  \, \left( \frac{2}{\pi k_0 r} \right)^\frac12 \exp (\iu k_0 r) \,
   \sum_{l=0}^{\infty}  \frac{g_l(\theta)}{r^l}
\end{equation}
Here $g_l(\theta)$ are some functions that absorb both $a_l$ and the $n$-index
summation and whose explicit form will not be needed.
The $\sim$ sign indicates that this well known series is,
as a more rigorous analysis shows, an asymptotic rather than necessarily
a convergent one \cite{Karp61,Zarmi-Turkel13}.

Even though functions $g_m$ depend on the solution and therefore are unknown \emph{a priori},
we still proceed and use the first few terms in \eqnref{eqn:u-eq-asymptotic-series}
as basis functions for our ``machine''. This works because $g_l$ depend only on the angle
$\theta$, while we deliberately choose the dof to be independent of $\theta$.
The general idea is best illustrated with a particular case of only two basis functions
$$
   \psi_0 = \frac{\exp(\iu k_0 r)}{\sqrt{k_0 r}} \, g_0(\theta),
   ~~~~~
   \psi_1 = \frac{\exp(\iu k_0 r)}{r \, \sqrt{k_0 r}} \, g_1(\theta)
$$
Since our dof need to be independent of $\theta$ (see above),
radial derivatives are a natural choice:
$$
   l_\beta (u) \,=\, \frac{\partial^\beta u}{\partial r^\beta},
   ~~~ \beta = 0,1,2
$$
Applying these dof to the basis set, one obtains by straightforward calculation
$$
   N^T \,=\, \{ l_\beta(\psi_\alpha )\}_{\alpha=1 ~ \beta=1}^{n ~~~~ m} \,=\,
   \frac{\exp(\iu k_0 r)}{\sqrt{k_0}} \,
   \begin{pmatrix}
   g_0(\theta) ~&~ 0\\
   0 ~&~ g_1(\theta)
   \end{pmatrix}
   \,
   \begin{pmatrix}
   -\frac{\iu}{r^{3/2}} ~&~
    \frac{2k_0 r + 3\iu}{2 r^{5/2}} ~&~
   -\frac{-\iu k_0^2 r^2 + 3k_0 r + 15\iu / 4}{r^{7/2}} \\
   \frac{1}{r^{1/2}}  ~&~
    \frac{-1 + 2\iu k_0 r}{2 r^{3/2}}  ~&~
   -\frac{k_0^2 r^2 + \iu k_0 r - 3/4}{r^{5/2}}
   \end{pmatrix}
$$
The null space of this matrix is seen to be independent of $\theta$,
and the coefficients for the absorbing condition are calculated to be
$$
  \underline{s} \,=\, \mathrm{null}~N^T \,=\,
  \left( \frac{3}{4 r^2} - k_0^2 - \frac{3 \iu k_0}{r}, ~~~
   \frac{3}{r} - 2\iu k_0, ~~~ 1 \right)^T
$$
More explicitly, the boundary condition is
$$
   \mathcal{L}u \,=\, 0,   ~~~~~
   \mathcal{L}u \equiv \sum_\beta s_\beta l_\beta(u) \,=\,
   \left( \frac{3}{4 r^2} - k_0^2 - \frac{3 \iu k_0}{r} \right) u \,+\,
   \left( \frac{3}{r} - 2\iu k_0 \right) \frac{\partial u}{\partial r}
   \,+\, \frac{\partial^2 u}{\partial r^2}
$$
which is none other than the second-order Bayliss-Turkel condition.

%%%%%%%%%%%%%%%%%%%%%%%%%%%%%%%%%%%%%%%%%
\section{Example: Sinusoidal Basis, Mixed Derivatives as dof}
\label{sec:Sin-basis-derivatives-dof}
%%%%%%%%%%%%%%%%%%%%%%%%%%%%%%%%%%%%%%%%%
%
Now consider the time-dependent wave equation \eqnref{eqn:del2-u-d2u-dt2}
in the half-plane $x > 0$. To run our ``machine,'' let us choose the basis of outgoing waves
\begin{equation}\label{eqn:psi-eq-dplane-dtheta}
   \psi_\alpha(x, y, t) = \frac{d^\alpha}{d \theta^\alpha} \,
   \exp \left(\iu (-x \cos \theta - y \sin \theta + t) \right)_{|\theta=0},
   ~~~ \alpha = 0,1, \ldots, n-1
\end{equation}
The rationale for this choice of functions is that they are expected to provide
accurate approximation of outgoing waves near normal incidence.
Explicit expressions for the first five of these functions are
$$
   \psi_0 = \exp(\iu (-x+t));
   ~~~
   \psi_1 = -\iu y \exp(\iu (-x+t));
   ~~~
   \psi_2 = (-y^2 + \iu x) \exp(\iu (-x+t))
$$
$$
   \psi_3 =  (\iu + \iu y^2 + 3x )y \, \exp(\iu (-x+t));
  ~~~~~
   \psi_4 = (-\iu x  - 3x^2 + 4 y^2  + y^4 - 6 \iu xy^2) \exp(\iu (-x+t))
$$
As dof, let us introduce
$$
  l_\beta(u) ~=~ d^{m_x}_x d^{m_y}_y d^{m_t}_t u
$$
where $m_x = 0, 1$; $m_y = 0, 2$; $m_x+m_y+m_t$ is either 2 (a second order method) 
or 3 (a third order method).
The omission of $m_y = 1$ reflects the symmetry of the problem with respect to $y$.

The respective $N^T$ matrices for $n=3$ and $n=5$ basis functions are
$$
   N^T \,=\, 
   \begin{pmatrix}
   -1 & 0 & 1\\
   0  &  0 &  0\\
   0  & -2 &  -1
   \end{pmatrix}
   ~~~~~~ (n=3)
$$
$$
   N^T \,=\, 
   \begin{pmatrix}
   -\iu & 0 & \iu & 0\\
   0  &  0 &  0 & 0\\
   0  & -2\iu &  -\iu & 2\iu\\
   0 &     0 &   0 &     0\\
   0 &   8\iu &  \iu &  -20\iu
   \end{pmatrix}
   ~~~~~~ (n=5)
$$
Calculating the null space of these matrices, one arrives at
the Engquist-Majda conditions of order two and three, respectively.
Thus not only the Bayliss-Turkel but also the Engquist-Majda conditions
can be generated by the proposed machine.

\textit{Remark}. Clearly, with an elementary degree of foresight,
the odd-numbered basis functions $\psi_{1,3}$ could have been omitted
from the basis set, as they produce zero rows of $N^T$ due to symmetry.
These functions were retained, however, to demonstrate the operation
of the Trefftz machine in a semi-automatic mode, with as little
``human intervention'' as possible.

%
%%%%%%%%%%%%%%%%%%%%%%%%%%%%%%%%%%%%%%%%%
\section{Sinusoidal Basis and Test Functions: a New Type of Boundary Condition}
\label{sec:Sin-basis-sin-test}
%%%%%%%%%%%%%%%%%%%%%%%%%%%%%%%%%%%%%%%%%
%
For the time-dependent problem \eqnref{eqn:del2-u-d2u-dt2},
with velocity $v$ normalized to unity,
let us now use a trigonometric basis of outgoing waves
$$
   \psi_{2\alpha}(x, y, t) ~=~ 
   \frac{d^\alpha}{d \theta^\alpha} \, \psi_{\cos}(x, y, t, \theta)_{| \theta = 0},  
   ~~~~
   \psi_{2\alpha-1}(x, y, t) ~=~ 
   \frac{d^\alpha}{d \theta^\alpha} \, \psi_{\sin}(x, y, t, \theta))_{| \theta = 0}
$$
where
$$
  \psi_{\cos}(x, y, t, \theta) \,=\, \cos(-x \cos \theta - y \sin \theta + t),
  ~~~ \psi_{\sin}(x, y, t, \theta) \,=\, \sin(-x \cos \theta - y \sin \theta + t)
$$
It is interesting to consider dof based on integrals rather than derivatives; for example:
\begin{equation}\label{eqn:lbeta-u-eq-int-int-int}
   l_\beta(u) ~=~ \int_{-\Deltarm}^{\Deltarm} \int_{-\Deltarm}^{\Deltarm} \int_{-\Deltarm}^{\Deltarm}
   u \, \psi_\beta \, dx \, dy \, dt
\end{equation}
Here, in the spirit of Galerkin methods, the test set $\psi$ coincides with the basis set.
$\Deltarm$ is an adjustable integration limit taken to be the same for all variables,
since for $v = 1$ the spatial and time scales are identical.

\begin{figure}
\centering
   \includegraphics[width=0.95\linewidth]{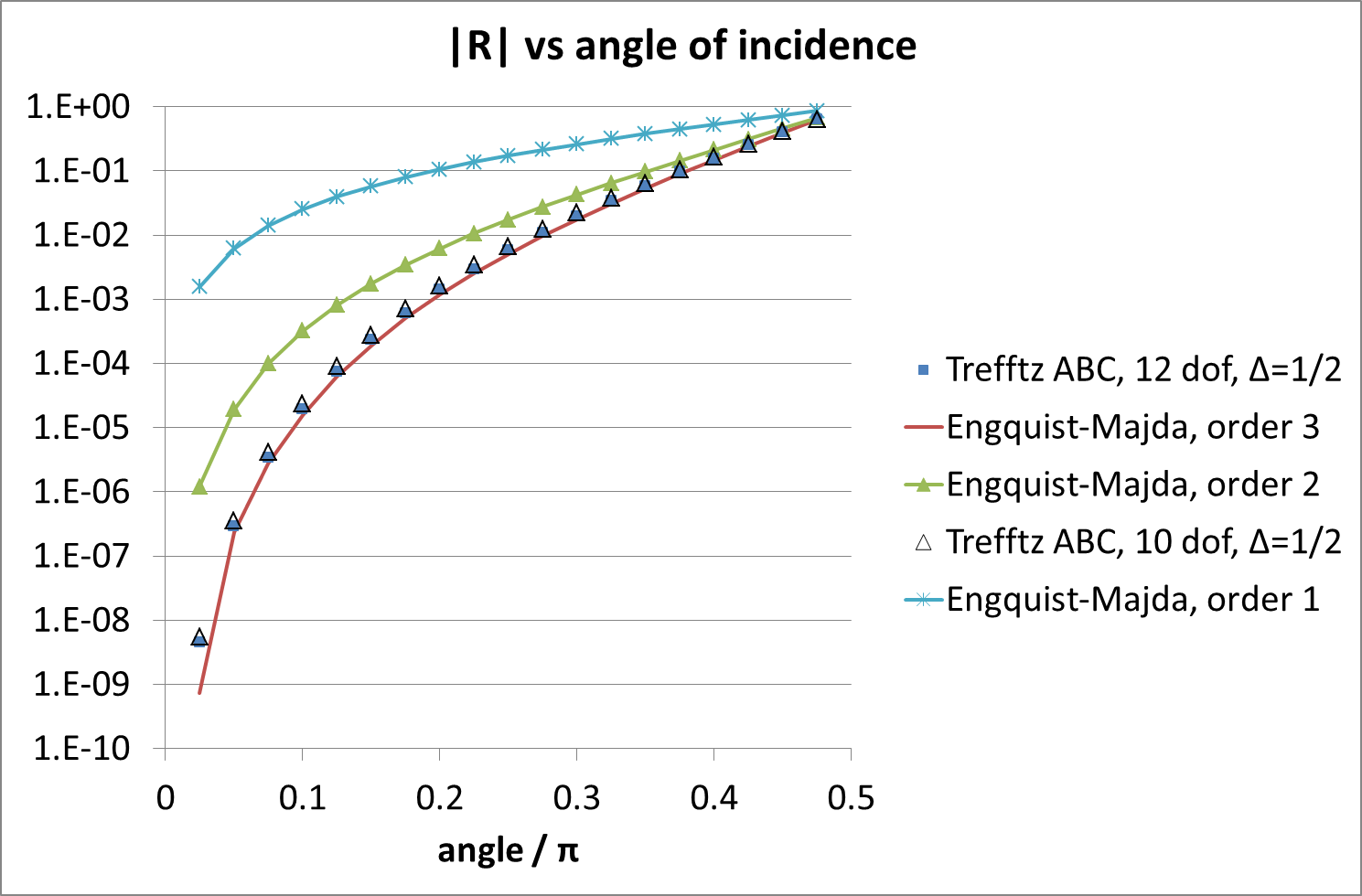}
    \caption{The absolute value of the reflection coefficient vs. the angle of incidence
    of a plane wave. Results for the ``trefftz machine'' with ainusoidal basis and test functions,
    in comparison with the Engquist-Majda conditions.}
    \label{fig:R-vs-angle-Trefftz-trig-basis-vs-EM}
\end{figure}

The absolute value of the reflection coefficient is plotted in \figref{fig:R-vs-angle-Trefftz-trig-basis-vs-EM}
as a function of the angle of incidence. One observes that the behavior of the
method with 10 or 12 integral dof \eqnref{eqn:lbeta-u-eq-int-int-int} 
is very close to that of the Engquist-Majda condition of order 3; 
however, the new method does not involve any derivatives.

%
%%%%%%%%%%%%%%%%%%%%%%%%%%%%%%%%%%%%%%%%%
\section{Three Absorbing Schemes in the Frequency Domain}
\label{sec:3schemes-FD}
%%%%%%%%%%%%%%%%%%%%%%%%%%%%%%%%%%%%%%%%%
%
In this section, we compare three absorbing schemes for the 2D Helmholtz 
equation \eqnref{eqn:del2-u-plus-k2u}
in the frequency domain. The domain is a square $\Omegarm = [-L, L] \times [-L, L]$.

\begin{enumerate}
  \item First, we consider the previously developed FLAME scheme \cite{Tsukerman06,Tsukerman-book07,Tsukerman05}
  over a six-point stencil on the sides of $\partial \Omegarm$ and over a four-point stencil at the corners
  of $\partial \Omegarm$. The basis set consists of five (on the straight sides) or three (for corner stencils)
  outgoing plane waves; see details below. The dof are, as in finite difference analysis,
  the nodal values of the solution.
  \item Same as above, but with the new basis set \eqnref{eqn:psi-eq-dplane-dtheta}. The rationale
  for this choice is to maximize the accuracy of approximation around normal incidence.
  However, approximation turns out to be good not only for very small angles
  but in a fairly broad range of angles of incidence.
  \item Integral dof. As basis functions $\psi_\alpha(x, y)$,
  we again use $\theta$-derivatives \eqnref{eqn:psi-eq-dplane-dtheta}
  of a plane wave at normal incidence, but the dof are now defined as integrals
  \begin{equation}\label{eqn:lu-eq-int-u-psi}
     l_\beta(u) \,=\, \int_0^{\Deltarm} u\, \psi_\alpha^* \, dx
  \end{equation}
  The test functions $\psi_\alpha$ in \eqnref{eqn:lu-eq-int-u-psi}
   are again the same as basis functions. Note that there is no $y$-integration in this example.
\end{enumerate}

The detailed setup of these three methods is as follows. In the first one (FLAME schemes),
the basis over the straight part $x=0$ of the boundary consists of five plane waves
$\psi_\alpha(x,y) = \exp(\iu k_0 (-x \cos \theta_\alpha - y \sin \theta_\alpha))$,
with $\theta_\alpha = \alpha \pi/6$, $\alpha = -2,-1,0,1,2$. The dof are the nodal values of these plane waves
on the six-point stencil $(x_\beta, y_\beta)$,  $\beta = 1,\ldots, 6$.
The coordinates of the stencil nodes are
$x_{1..6} = \{0, h, 0, h, 0, h \}$, $y_{1..6} = \{0, 0, -h, -h, h, h\}$, where for simplicity
the origin is set at the first node and $h$ is the grid size.
The coefficient vector of the FLAME scheme, i.e. of the absorbing condition,
is $\underline{s} = \mathrm{Null}~ N^T$, where $N^T_{\alpha \beta} = \psi_\alpha(x_\beta, y_\beta)$.
Expressions for these coefficients are too cumbersome to be listed here
but easily obtainable with symbolic algebra.   For reference, the numerical values of
these coefficients for $h = \lambda_0/12 = 2\pi/(12 k_0)$ are
%
%\begin{center}
%\begin{tabular}{c}
  $s_1 = 0.35149777 - 1.3321721 \iu$,
  $s_2 = -1.36164086 - 0.21016073 \iu$,
  $s_3 = -0.39962632 + 0.91667814 \iu$,
  $s_4 =  1$,
  $s_5 = -0.39962632 + 0.91667814 \iu$,
  $s_6 = 1$.
%\end{tabular}
%\end{center}

At a corner (placed for simplicity at the origin), the four-point stencil
is $x_{1..4} = \{0, h, 0, h\}$, $y_{1..4} = \{0, 0, h, h\}$, and the three basis functions
are $\psi_\alpha(x,y) = \exp(\iu k_0 (-x \cos \theta_\alpha - y \sin \theta_\alpha))$,
$\theta_\alpha = \alpha \pi/6$, $\alpha = 0,1,2$.

In Method 2 above, the grid stencils are the same as in Method 1,
but with the basis set \eqnref{eqn:psi-eq-dplane-dtheta}. The absorbing scheme
is again found as the null space of the respective matrix $N^T$, although
this matrix is of course different from that of Method 1. The scheme
is simple enough to be written out explicitly:
$$
   s_1 = -2 \exp(-\iu k_0 h) \, \frac{5(k_0 h)^2 + 3\iu k_0 h - 3}
   {(k_0 h)^2 + 3\iu k_0 h + 3};
~~~
   s_2 =  -2 \, \frac{-5(k_0 h)^2 + 3\iu k_0 h + 3}
   {(k_0 h)^2 + 3\iu k_0 h + 3};
$$ 
$$
   s_3 = s_5 = -\exp(-\iu k_0 h) \, \frac{(k_0 h)^2 - 3\iu k_0 h + 3}
   {(k_0 h)^2 + 3\iu k_0 h + 3};
~~~
   s_4 = s_6 = 1
$$

The respective scheme at the corner for Method 2 is
$$
  s_{1,\mathrm{corner}} = -s_{3,\mathrm{corner}} =  \exp(-\iu k_0 h);  ~~
  s_{4,\mathrm{corner}} = -s_{2,\mathrm{corner}} =  1
$$

In Method 3, the basis functions $\psi_\alpha(x, y)$
are defined as $\theta$-derivatives \eqnref{eqn:psi-eq-dplane-dtheta}
of a plane wave at normal incidence, but the dof are now defined as integrals
\eqnref{eqn:lu-eq-int-u-psi}.

\begin{figure}
\centering
   \includegraphics[width=0.95\linewidth]{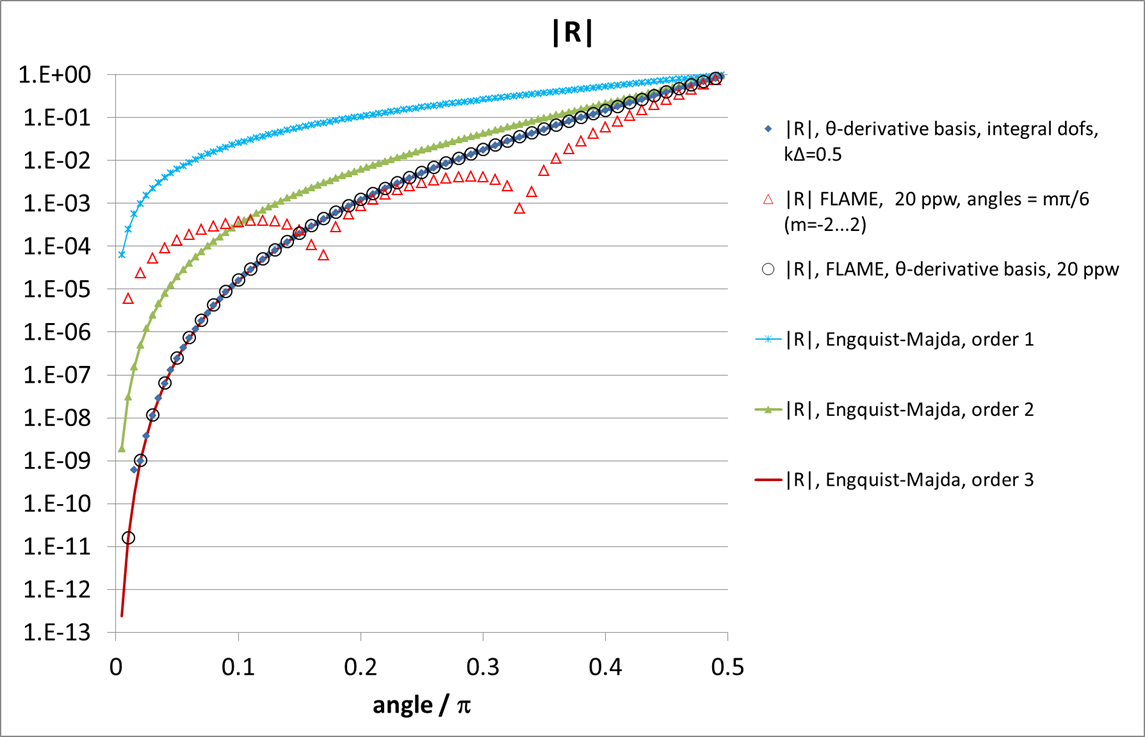}
    \caption{The absolute value of the reflection coefficient vs. the angle of incidence
     of a plane wave. Results for the ``Trefftz machine'' with ainusoidal and
     $\theta$-derivative bases, in comparison with the Engquist-Majda conditions.}
    \label{fig:R-vs-angle-Trefftz-FLAME-vs-EM}
\end{figure}

The absolute value of the reflection coefficient for all three methods is plotted
in \figref{fig:R-vs-angle-Trefftz-FLAME-vs-EM} (20 points per wavelength,
i.e. $\lambda_0 / h = 20 \Leftrightarrow k_0 h = \pi/10$); the respective results for
the Engquist-Majda conditions of orders one through three are also shown for reference.
It is evident that $|R(\theta)|$ in Methods 2 and 3, which are both based on the $\theta$-derivative
basis \eqnref{eqn:psi-eq-dplane-dtheta}, are virtually indistinguishable from
that of the Engquist-Majda condition of order three.
Admittedly, Method 3 is an analytical integral condition
whose performance will degrade once the integral is approximated in conjunction
with a given discretization scheme (finite difference or finite element).
At the time of this writing, I am not aware of a discrete version of Method~3
whose performance would be on a par with that of the analytical condition.
For this reason, we shall not discuss Method 3 further
and turn now to a comparison of Methods 1 and 2.

Since basis functions in Method 2 are tailored toward approximation near normal incidence,
it is not surprising that this method outperforms Method 1 for small angles of incidence
($\theta \lesssim \pi /6$, \figref{fig:R-vs-angle-Trefftz-FLAME-vs-EM}).
At greater angles, it is Method 1 that yields lower reflection. In practical simulations,
one may therefore expect that if the artificial boundary is placed far away
from the scatterers and consequently the scattered field impinges on it
at close-to-normal incidence, Method 2 will be preferable; otherwise Method 1
can be expected to perform better.

% % % % % % % % % % % % % % % % % % % % % % % %
\section{Conclusion}
\label{sec:Conclusion}
% % % % % % % % % % % % % % % % % % % % % % % %
%
The paper develops a generator of high-order nonreflecting boundary conditions for
wave problems. This generator is based on a set of local Trefftz basis functions
(outgoing waves) and a commensurate set of linear functionals (degrees of freedom).
Degrees of freedom involving nodal values on a grid give rise to numerical
(finite-difference-type) nonreflecting conditions, while dof involving derivatives or integrals
produce analytical ones. The schemes, analytical as well as numerical,
are given by the simple nullspace formula \eqnref{eqn:s-in-null-Nt}.
Consistency of such schemes can be established in a way similar to the analysis
of \cite{Tsukerman06,Tsukerman-book07}; stability and convergence cannot be guaranteed
\textit{a priori} and need to be examined on a case-by-case basis.
Nevertheless classical boundary conditions such as Engquist-Majda and Bayliss-Turkel,
and likely also their extensions \cite{Zarmi-Turkel13},
can be reproduced faithfully by the proposed ``Trefftz machine''.
Corners and edges are treated algorithmically
the same way as straight boundaries. The proposed approach opens up
various avenues for the development of new approximate boundary conditions
and gives an opportunity to look at the existing ones from a different perspective.
Extensions to 3D problems, problems with frequency dispersion, and to Maxwell's electrodynamics
are certainly possible.

%%%%%%%%%%%%%%%%%%%%%%%%%%%%%%%%%%%%%%
\section*{Acknowledgment}
\label{sec:Acknowledgment}
%%%%%%%%%%%%%%%%%%%%%%%%%%%%%%%%%%%%%%
%
I thank Professor Dmitry Golovaty for discussions.
\bibliographystyle{unsrt}
\bibliography{Igor_reference_dbase}
\end{document}